\documentclass[12pt]{amsart}

\usepackage{amsfonts,amsmath,latexsym,amssymb,verbatim,amsbsy,times}
\usepackage{amsthm}
\usepackage{pstricks}
\theoremstyle{plain}
\newtheorem{THEOREM}{Theorem}[section]

\newtheorem{theorem}[THEOREM]{Theorem}
\newtheorem{corollary}[THEOREM]{Corollary}
\newtheorem{lemma}[THEOREM]{Lemma}

\theoremstyle{definition}

\newtheorem{definition}[THEOREM]{Definition}

\theoremstyle{remark}

\newcommand{\thm}[1]{Theorem~\ref{#1}}

\newcommand{\R}{\ensuremath{\mathbb{R}}}   %%% reals
\def \a {\alpha}
\def \b {\beta}
\def \d {\delta}
\def \g {\gamma}
\def \G {\Gamma}
\def \e {\epsilon}
\def \f {\varphi}

\def \n {\nabla}
\def \s {\sigma}

\def \t {\tau}

\def \O {\Omega}

\def \reg {\mathcal{R}}
\def \D {\Delta}

\def \< {\langle}
\def \> {\rangle}
\def \p {\partial}
\def \ra {\rightarrow}
\def \ss {\subset}

\newcommand{\sgn}[1]{\mathrm{sgn}(#1)}

\def \loc {\mathrm{loc}}

\DeclareMathOperator{\supp}{supp} %
\DeclareMathOperator{\diver}{div} %
\DeclareMathOperator{\tr}{Tr} %

\begin{document}

\title{On the energy of inviscid singular flows}
\author{Roman Shvydkoy}
\thanks{The work was partially supported by NSF grant DMS -- 0604050. The author is grateful to R. Caflisch and A. Cheskidov for stimulating discussions.}
\address{Department of Mathematics, Stat. and Comp. Sci.\\
 851 S  Morgan St., M/C 249\\
        University of Illinois\\
        Chicago, IL 60607}
\email{shvydkoy@math.uic.edu}

\date{\today}

\begin{abstract}
It is known that the energy of a weak solution to the Euler equation is conserved if it is slightly more regular than the Besov space $B^{1/3}_{3,\infty}$.
When the singular set of the solution is (or belongs to) a smooth manifold, we derive various $L^p$-space regularity criteria dimensionally equivalent to the critical one. In particular, if the singular set is a hypersurface the energy of $u$ is conserved provided the one sided non-tangential limits to the surface exist and the non-tangential maximal function is $L^3$ integrable, while the maximal function of the pressure is $L^{3/2}$ integrable. The results directly apply to prove energy conservation of the classical vortex sheets in both 2D and 3D at least in those cases where the energy is finite.
\end{abstract}



\maketitle

\section{Introduction}

In this paper we study weak solutions to the Euler equations modeling evolution of inviscid fluid
flows
\begin{align} \label{e:Euler}
\frac{\partial u}{\partial t} + (u \cdot \nabla)u &= - \nabla p, \\
\nabla \cdot u &=0. \label{diver-free}
\end{align}
Here $u$ is a divergence-free velocity field, and $p$ is the internal pressure. The classical law of energy conservation
$$
\int |u(t)|^2 dx = \int |u_0|^2 dx
$$
for smooth rapidly decaying solutions of \eqref{e:Euler} and \eqref{diver-free} is an easy consequence
of the antisymmetry of the nonlinear term. Weak solutions to \eqref{e:Euler}
are believed to describe turbulent phenomena at large Renolds number in the inertial range of frequencies. The Kolmogorov-Obukhov power laws predict
solutions to be $\frac{1}{3}$-H\"{o}lder continuous in a statistically averaged sense. Moreover, since the energy is not lost within the inertial range the energy
flux through inertial scales is to be proportional to the mean energy dissipation rate $\e$ (\cite{Frisch}). Experiments show that $\e$
is essentially independent of the viscosity coefficient. So, if in the limit of infinite Reynolds number turbulent solutions converge
in some sense to weak solutions of the Euler solutions, then such solutions are expected to be on average energy dissipative.

Onsager \cite{Onsager} stated that all $(\frac{1}{3}+\d)$-regular solutions conserve energy, and there may exist solutions exactly $\frac{1}{3}$-regular that do not.
The results of Eyink \cite{Eyink} followed by the work of Constantin, E and Titi \cite{CET} give Onsager's hypothesis rigorous proof in the spaces $B^{1/3+\d}_{3,\infty}$, which measure H\"{o}lder continuity in the $L^3$-space. An example of a vector field exhibited in \cite{EggersGrossmann,Eyink} suggests that the exponent $\frac{1}{3}$ may indeed be critical, however no rigorous proof of this fact exists at the moment. An improvement upon \cite{CET} by Duchon and Robert \cite{Duchon} showed that some solutions conserve energy even in the Onsager-critical case. In recent paper \cite{ccfs} the criterion was established in the dimensionally optimal regularity class $L^3_t B^{1/3}_{3,c_0}$ where $c_0$ signifies the decay $2^q \|\D_q u \|_3^3 \ra 0$ of the $\frac{1}{3}$-derivatives of the dyadic parts.

This present paper is motivated by the work of Caflisch, Klapper and Steele \cite{CKS}, where the authors obtain bi-H\"{o}lder sufficient conditions for solutions with singularity set located on a smooth submanifold of $\R^n$. Although these conditions are subcritical, they are more practical in applications, for example, to multifractal models of turbulence  (see \cite{CKS,FrischP}). However, other important classes of singular weak solutions such as vortex sheets remain unattainable by the results of \cite{CKS,ccfs,CET}. Indeed, classical analytic vortex sheets in 2D or in 3D fall exactly into the critical class $B^{1/3}_{3,\infty} \backslash B^{1/3}_{3,c_0}$.

In this paper we study the energy law for solutions which exhibit organized singular sets. Examples of singular set organization include time dependent families of submanifolds of $\R^n$ and their locally finite unions. We obtain Onsager-critical criteria near such sets in terms of $L^p$-spaces, which do not involve calculation of spacial H\"{o}lder exponents. For instance, in the case of a 3D solution with point singularity $s(t)$ at  time  $t\in [0,T]$ and $s \in C^{3/5}([0,T])$ the energy of $u$ is conserved provided $u \in L^3_t L^{9/2}_x$ near the curve $s$ (see also application to viscous flows in \cite{s}). In higher dimension we use mixed $L^p$-spaces relative to the singular manifold (see \thm{T:main} and Section \ref{SS:oo}). The case of hypersurface $S(t)$ is treated separately in Section \ref{S:slits}. We will introduce the notion of a slit suitable for subsequent analysis. We assume that the velocity and pressure fields have non-tangential or normal limits and that the non-tangential maximal functions are integrable on the surface. As a consequence of weak formulation of the Euler equations, we show that  all  slits necessarily satisfy the kinematic condition similar to that of a free surface, so that particles that are initially  on the surface  stay on the surface at all time (see Lemma \ref{L:slit}). This case is radically different from the lower dimensional case where no particular evolution law is imposed by the equation. Our analysis shows that the energy of a solution $u$ with a slit type of singularity is conserved provided the non-tangential maximal functions of $u$ and the pressure $p$ belong to $L^3(S)$ and $L^{3/2}(S)$, respectively (see \thm{T:slit}). These conditions are verified for the classical 2D and 3D vortex sheets in Section \ref{S:sheets} implying their energy conservation (under zero total circulation in 2D).

Energy non-conservative weak solutions without any apparently organized space singularities have long been constructed by Scheffer \cite{Scheffer} and Shnirelman \cite{Shnirelman}, and more recently by De Lellis and Sz\'{e}kelyhidi in \cite{DeLellis}. Those belong to $L^2_tL^2_x$  and $L^\infty_tL^2_x$, respectively, and therefore are considerably Onsager-supercritical. As we mentioned earlier the vector field considered by Eyink \cite{Eyink} with non-vanishing energy flux belongs exactly to $B^{1/3}_{3,\infty} \backslash B^{1/3}_{3,c_0}$. However, no weak solution with this initial condition is known to exist. The example serves to show that the traditional mollification argument used to prove energy conservation is sharp. Again, one can show that it has no organized singularities. It is in fact locally nowhere in the energy-regular class $B^{1/3}_{3,c_0}$.

Although we chose to use $\R^n$ as a model case, the local nature of the arguments presented below allows us to
apply the results to other boundary problems, such as periodic in all or some spacial directions. This will be especially useful in application to vortex sheets.

\section{Weak solutions and Regular sets}\label{S:reg}

\begin{definition}\label{D:weaksol}
A vector field $u\in C_w([0,T]; L^2(\R^n))$, (the space of weakly continuous functions), is a weak solution of the Euler equations
with initial data $u_0\in L^2(\R^n)$
if for every $\psi \in C^\infty_0([0,T]\times \R^n)$
with $\n_x \cdot \psi =0$ and $0\leq t \leq T$, we have
\begin{equation}\label{weaksol}
    \int_{\R^n \times \{t\}} u \cdot \psi - \int_{\R^n \times \{0\}} u_0 \cdot \psi - \int_0^t \int_{\R^n} u\cdot
    \p_s \psi = \int_0^t \int_{\R^n} (u\otimes u) : \n \psi,
\end{equation}
and $\n_x \cdot u(t) = 0$ in the sense of distributions. We define the operation $:$ by
$$
A:B = \tr[AB].
$$
\end{definition}

It will be convenient to work with the associated pressure defined by
\begin{equation}\label{pressure}
p = - \sum_{l,k=1}^n R_lR_k(u_lu_k),
\end{equation}
where $R_l$ are the classical Riesz projections. With the use of $p$ we can alternatively restate the definition of a
weak solution without requiring $\n_x \cdot \psi =0$. Namely,
\begin{multline}\label{weaksolalt}
    \int_{\R^n \times \{t\}} u \cdot \psi - \int_{\R^n \times \{0\}} u \cdot \psi - \int_0^t \int_{\R^n} u\cdot
    \p_s \psi \\= \int_0^t \int_{\R^n} (u\otimes u) : \n \psi + \int_0^t \int_{\R^n} p \diver \psi,
\end{multline}
holds for all $\psi \in C^\infty_0([0,T]\times \R^n)$. Since the pressure is only a distribution, the pairing between $p$ and $\diver \psi$ is to be understood
accordingly.

Based on the results of \cite{ccfs} we introduce the global regularity class $\reg(\R^n \times I)$ consisting of vector fields $u \in L^3(\R^n \times I)$ on a time interval $I\ss [0,T]$ such that
\begin{equation}\label{reg}
\lim_{y \ra 0} \frac{1}{|y|} \int_{\R^n\times I} |u(x-y,t) - u(x,t)|^3 dx dt = 0.
\end{equation}
 For an open set $U\ss \R^n$ is an open set, we define $\reg(U \times I)$ as the class of fields $u$ such that $u \phi \in \reg(\R^n \times I)$ for all $\phi \in C_0^\infty(U)$.

Alternatively, we could define $\reg(\R^n \times I)$  using Littlewood-Paley decomposition over dyadic shells in the frequency space (see \cite{Stein})
$$
u = \sum_{q=0}^\infty \D_q u.
$$
Thus, condition \eqref{reg} is equivalent to
\begin{equation}\label{reg:equiv}
\lim_{q \ra \infty} \int_I 2^q \|\D_q u(t)\|_3^3 dt = 0.
\end{equation}
In this form the regularity class was introduced in \cite{ccfs}, and the energy conservation was established. A similar but less time-optimal
class was considered in \cite{Duchon} as a direct improvement upon \cite{CET}. We remark that condition $u\in L^3([0,T]; B^{1/3}_{3,c_0})$ implies \eqref{reg:equiv}, where $c_0$ stands to indicate
$$
\lim_{q \ra \infty} 2^q \|\D_q u(t)\|_3^3 = 0.
$$

\begin{definition}\label{D:reg} Let $u$ be a weak solution to the Euler equations.
A point $(x_0,t_0) \in \R^n\times [0,T]$ is called regular if there exists  an open neighborhood $U\ss \R^n$ of $x_0$
and a relatively open interval $I \ss [0,T]$ containing $t_0$ such that $u \in \reg(U\times I)$.
An open set $D \ss \R^n \times [0,T]$ is regular if every point in it is regular. The set $S$ of all irregular points  is called the singular set
of $u$.
\end{definition}

The main purpose of this section is to prove the following local energy balance relation
inside every regular set. For a set $A \ss \R^n \times [0,T]$ we denote by $A(t)$ the slice
$A \cap \R^n \times \{t\}$.

\begin{lemma}\label{L:reg}
Let $D$ be a regular set of a weak solution $u$. Then for every $\phi \in C_0^\infty(D)$
one has
\begin{equation}\label{E:reg}
\int_{D(t'')} |u|^2 \phi - \int_{D(t')}|u|^2 \phi - \int_D |u|^2 \p_t \phi = \int_D (|u|^2 + 2p) u\cdot \n \phi,
\end{equation}
for all $t',t'' \in [0,T]$.
\end{lemma}

Before we prove this lemma, we need to take another seemingly obvious but not entirely straightforward step
by showing that one can substitute a mollified in space solution $u$ into \eqref{weaksol} as a test function. This fact is not so straightforward
since a priori $u$ may not have sufficient time regularity. The difficulty has been removed in a similar situation in \cite{Kukavica} by
considering mollification both in space and time, however in our case such mollification would introduce unnecessary technical obscurity.
So, let us fix a mollifier $h \in C_0^\infty(\R^n)$ with $\int h = 1$ and $h = 0$ outside the unit ball. Denote
\begin{align*}
u_\d(x,t) & = \int_{\R^n} h_\d(y) u(x-y,t) dy,\\
h_\d(y) & = \d^{-n} h(y \d^{-1}).
\end{align*}
We prove the following lemma.

\begin{lemma}\label{L:moll}
Let $u$ be a weak solution. Then for each fixed $\d>0$, $u_\d:[0,T] \ra W^{s,q}$ is absolutely continuous for all $s \geq 0$ and $q\geq
2$, and moreover
\begin{equation}\label{E:moll}
\p_t u_\d  = - \n \cdot (u \otimes u)_\d - \n p_\d,
\end{equation}
for a.e. $t \in [0,T]$.
\end{lemma}
\begin{proof}
Substituting test-functions of the form
$$
\psi = \b(t) \psi(x),
$$
where $\b \in C_0^\infty(0,T)$ and $\psi \in C_0^\infty(\R^n)$ into \eqref{weaksolalt} we obtain
$$
\p_t \int u(t) \cdot \psi = \int u(t) \otimes u(t) :  \n \psi + \int p(t) \diver \psi
$$
in the distributional sense.
Hence, since $u$ is weakly continuous
$$
\int u(t) \cdot \psi = \int u(0) \cdot \psi + \int_0^t \int u(s) \otimes u(s) : \n \psi ds + \int_0^t \int p(s) \diver \psi ds,
$$
for all $0\leq t \leq T$. Let $\psi_l$ denote the coordinate components of $\psi$.
Taking the Fourier transform we obtain
\begin{multline*}
\int \hat{u}(t) \cdot \overline{\hat{\psi}} = \int \hat{u}(0) \cdot
\overline{\hat{\psi}} - i \int_0^t \int_{\R^n} (u_l
u_k)^\wedge(\xi,s)\xi_k  \overline{\hat{\psi_l}}(\xi) d\xi ds \\
- i \int_0^t \int_{\R^n} \hat{p}(\xi,s) \xi_l \overline{\hat{\psi_l}}(\xi) d\xi ds ,
\end{multline*}
assuming the usual summation convention. Let us notice that $(u_l
u_k)^\wedge$ and $\hat{p}$ are continuous and bounded functions of $\xi\neq 0$ for every $s$.
Let $\Sigma_t$ denote  the common Lebesgue set
of $\hat{u}(t)$ and $\hat{u}(0)$ not containing the origin, so that $| \R^n \backslash \Sigma_t | = 0$. Denote by
$e_j(\xi)$, $j =
1,...,n$ the vectors of the standard unit basis. For every $j$ and $\xi \in \Sigma_t$ we apply the previous identity to a sequence of functions $\psi$
such that  $ \overline{\hat{\psi^n}}(\xi) \ra e_j(\xi) \d_0(\cdot - \xi)$, where $\d_0$ is the Dirac mass.
We obtain
\begin{equation*}\label{alpha-version}
\hat{u}_j(\xi,t) \cdot e_\a(\xi) = \hat{u}_j(\xi,0) \cdot e_\a(\xi) - i  \int_0^t  (u_j u_k)^\wedge(s,\xi)\xi_k ds - i  \int_0^t \hat{p}(\xi,s) \xi_j ds,
\end{equation*}
for all $t \in [0,T]$ and $\xi \in \Sigma_t$. Thus, the identity
\begin{equation}\label{distr1}
u(t) = u_0 - \int_0^t [ \n \cdot (u \otimes u) + \n p ] ds
\end{equation}
holds in the sense of distributions for all $t\in [0,T]$.
Mollifying \eqref{distr1} with $h_\d$ we obtain
\begin{equation}\label{distr2}
u_\d(t) = u_\d(0) - \int_0^t [ \n \cdot (u \otimes u)_\d + \n p_\d ] ds
\end{equation}
for all $t \in [0,T]$. Since $u(0) \in L^2$ we have $u_\d(0) \in
W^{s,q}$ for all $s \geq 0$ and $q \geq 2$, and since $u\otimes u
\in  L^\infty_tL^1_x$, we have $ \n \cdot (u \otimes u)_\d,\ \n p_\d \in
L^\infty_t W^{s,q}_x$, for all $s \geq 0$ and $q \geq 1$.
This proves the lemma.
\end{proof}

Let us denote
$$
B_c = (-c,c)^{n}.
$$

\begin{proof}[Proof of Lemma \ref{L:reg}]

First let us observe that $p\in L^{3/2}_\loc(D)$. Indeed, for a compact subset $K \ss D$ let $\e>0$ be such that
$\overline{K+B_\e} \ss D$. Let $\a \in C_0^\infty(D)$ be such that $\a \equiv 1$ on $\overline{K+B_{\e/2}}$ and
$\a \equiv 0$ on $D \backslash K+ B_\e$. Then
\begin{equation}\label{pressuredec}
p = R_iR_j(u_iu_j \psi) +R_iR_j(u_iu_j (1-\a)).
\end{equation}
Since $ u_i u_j \a \in L^{3/2}$, so is the first term in \eqref{pressuredec}. The second term
belongs to $L^\infty(K)$ since $1-\a \equiv 0$ on $K+B_{\e/2}$, $u_iu_j \in L^1$ and the kernel of $R_iR_j$ is bounded away from the $\e/2$-neighborhood of the origin. This observation justifies the pressure integral in \eqref{E:reg}.

Using partition of unity over the support of $\phi$ we reduce the lemma to the case
$D = U\times I$, where $U$ is an open ball. So, suppose $\phi \in C^\infty_0(U\times I)$. Choose $\d_0>0$ so small
that $$\overline{\supp(\phi(\cdot,t) + B_{\d_0})} \ss K \ss U,$$ for all $t\in I$. Let us now use \eqref{E:moll} with $\d < \d_0$.
We obtain
$$
\p_t u_\d \cdot u_\d \phi = \frac{1}{2} \p_t(|u_\d|^2 \phi) -
\frac{1}{2}|u_\d|^2 \p_t \phi.
$$
Integrating in time on $[t',t'']\ss [0,T]$ we obtain
\begin{multline}\label{distr4}
 \int_{U \times \{t''\}} |u_\d|^2 \phi -  \int_{U \times \{t'\}} |u_\d|^2 \phi - \int_{U \times[t',t'']} |u_\d|^2 \p_t
\phi
 = \\ = 2 \int_{U \times [t',t'']} (u\otimes u)_\d : \n(u_\d \phi) + 2 \int_{U \times [t',t'']} p_\d u_\d \cdot \n \phi .
\end{multline}
Notice that the time integration is in fact happening on the interval $[t',t''] \cap I$. So, we can pass to the limit as $\d \ra 0$ on the left hand side and in the pressure term. The nonlinear term will be treated similar to
\cite{CET}. First, consider a scalar $\b \in C_0^\infty(U)$ with $\b \equiv 1$ on $K$. We can then replace $u$ by $u\b$ under the integrals of \eqref{distr4}. Without
further change of notation we simply assume $u \in \reg(\R^n \times I)$.
We have
\begin{multline*}
2 \int_{U \times [t',t'']} (u\otimes u)_\d : \n(u_\d \phi) = 2 \int_{U \times [t',t'']} (u\otimes u)_\d : (u_\d \otimes \n \phi)\\ + 2 \int_{U \times [t',t'']} (u\otimes u)_\d : \n(u_\d) \phi.
\end{multline*}
Clearly, we can pass to the limit
\begin{equation}\label{addup}
2 \int_{U \times [t',t'']} (u\otimes u)_\d : (u_\d \otimes \n \phi) \ra 2 \int_{U \times [t',t'']} |u|^2 u \cdot \n \phi.
\end{equation}
Let us observe the following identity
\begin{equation}
(u\otimes u)_\d = r_\d(u,u)  - (u - u_\d) \otimes (u - u_\d) +  u_\d \otimes u_\d,
\end{equation}
where
$$
r_\d(u,u)(x,t) = \int_{\R^n} h_\d(y) (u(x-y,t)-u(x,t)) \otimes  (u(x-y,t)-u(x,t)) dy.
$$
Notice
$$
u(x) - u_\d(x) = \int h_\d(y) (u(x) - u(x-y)) dy,
$$
and
$$
\n u_\d (x) =  \frac{1}{\d} \int (\n h)_\d(y) \otimes (u(x) - u(x-y)) dy.
$$
So, we can estimate using H\"{o}lder and Minkowski inequalities
\begin{align*}
&\left| \int (u - u_\d) \otimes (u - u_\d) : \n(u_\d) \phi \right| \\& \leq \left( \int_{\R^n} h_\d(y) \int_{U\times I} |u(x,t) - u(x-y,t)|^3 dx dt dy \right)^{2/3}
\\ & \times \left( \int_{\R^n} \frac{1}{\d} |(\n h)_\d(y)| \int_{U\times I} |u(x,t) - u(x-y,t)|^3 dx dt dy \right)^{1/3}
\\ & \leq \frac{o(\d)}{\d^{1/3}} \left( \int |y|h_\d(y) dy \right)^{2/3} \left( \int |y| |(\n h)_\d(y)| dy \right)^{1/3} \leq o(\d) \ra 0.
\end{align*}
Similarly, the term with $r_\d$ vanishes as well. Finally,
\begin{multline*}
2 \int u_\d \otimes u_\d : \n (u_\d ) \phi = - \int |u_\d|^2 u_\d \cdot \n \phi   \ra  - \int |u|^2 u \cdot \n \phi.
\end{multline*}
This adds up with \eqref{addup} to produce the corresponding term in \eqref{E:reg}.
\end{proof}

\section{Low-dimensional singular sets}\label{S:low}

\begin{definition}\label{D:cover}
We say that a set $S \ss \R^n\times [0,T]$ admits
a \emph{$k$-dimensional $C^{\g,1}$-cover} if for every point
$(x_0,t_0)$ in the space-time there is an open neighborhood $U$ of $x_0$
in $\R^n$ and a relatively open subinterval
$I \ss [0,T]$ containing $t_0$ for which there exists a family of $C^1$-diffeomorphisms
\begin{equation}\label{diffeo}
\f_t: U \ra B_1, \quad t \in I,
\end{equation}
 satisfying the following conditions
 \begin{itemize}
 \item[(a)] $S(t) \cap U \ss \f_t^{-1}( \R^k \times \{0\}^{n-k} \cap B_1)$, for all $t\in I$;
 \item[(b)] There is $C>0$ such that
 \begin{equation*}
 \sup_{x\in U} |\f_{t'}(x) - \f_{t''}(x)| \leq C|t' - t''|^\g,
 \end{equation*}
 for all $t',t''\in I$;
 \item[(c)]  $\sup_{x \in U, t\in I} |\n_x \f_t(x)| \leq C$.
\end{itemize}
\end{definition}

\begin{theorem}\label{T:main}
Let $u \in L^{3}(\R^n \times [0,T])$ be a weak solution to the Euler equation on the time
interval $[0,T]$. Then $u$ conserves energy provided the singular set $S$ of $u$ admits a $k$-dimensional
$C^{\g,1}$-cover and  $u\in L^3 L^q_\loc$, where the values of $\g, n, k, q > 0$ satisfy
\begin{equation}\label{E:expon}
\g \geq \frac{q}{(q-2)(n-k)},\ n\geq k+2,\ q \geq 3\frac{n-k}{n-k-1}.
\end{equation}
\end{theorem}
\begin{proof}
We claim that in order  to prove \thm{T:main} it suffices to show
that for every coordinate chart $ U \times I$ and scalar
test-function $\phi \in C_0^\infty(U)$ independent of time one has
the following identity
\begin{equation}\label{localee}
 \int_{U \times \{t''\}} |u|^2 \phi -  \int_{U \times \{t'\}} |u|^2 \phi = \int_{U \times [t',t'']} ( |u|^2 +2 p ) u \cdot \n
 \phi,
\end{equation}
for all $t',t'' \in I$. Indeed, if this is the case, we fix an
arbitrary smooth $\phi$ with $\supp(\phi) \ss B_R$, and $t_0 \in
[0,T]$. By compactness we can find a finite collection of charts
$U_i \times I_i$, $i=\overline{1,M}$ so that all $I_i$'s contain
$t_0$ and $U_i$'s cover $B_R$. Put $I_0 = \cap_{i=1}^M I_i$.
Consider a partition of unity $\{g_i\}_{i=1}^M$ subordinate to the
cover, so that $\supp g_i \ss U_i$ and $\sum_{i=1}^M g_i = 1$ on
$B_R$. Since we have \eqref{localee} for any $\phi g_i$ and $t',t''
\in I_0$ summing up over $i$ we obtain \eqref{localee} for the given
$\phi$ itself. The above construction is carried out for every $t_0
\in [0,T]$. Thus, we can find a finite cover of $[0,T]$ by intervals
such as $I_0$, and as a consequence obtain \eqref{localee} for all
$t',t'' \in [0,T]$. Letting $\phi = \phi_0(x/R)$, where $ \phi_0 =
1$ on $B_1$ and $\phi_0 =0$ on $B_2$, and letting $R\ra \infty$ we
see that the right hand side of \eqref{localee} vanishes and we
arrive at the desired energy equality.

We will prove \eqref{localee} with the
use of Lemma \ref{L:reg}, but first we need to introduce a cut-off of
the singular sets $S(t)\cap U$. Let $\f_t : U \ra B_1$ be the
coordinate map, for $t\in I$. Denote $\overline{I} = [a,b]$. If
$\f_t$ is not defined at $a$ or $b$, then $t_0$ is not that point.
In this case we can consider a slightly shorter interval $I$ still
containing $t_0$ and so that $\f_t$ is defined at both ends. Let us
define an extension of $\f_t$ as follows
\begin{equation}\label{ext}
\tilde{\f}_t = \left\{%
\begin{array}{ll}
    \f_a, & t \leq a \\
    \f_t, & a < t < b \\
    \f_b, & t \geq b\\
\end{array}%
\right.
\end{equation}
Notice that $\tilde{\f}_t$ still satisfies condition (b) of
Definition \ref{D:cover} on the entire real line. Let $\b(\t)$ be a
mollifier. Define
$$
\f_{t,\e}(x) = \int_\R \e^{-1} \b(\t \e^{-1})\tilde{\f}_{t-\t}(x)
d\t.
$$
Let us notice the following approximation inequalities:
\begin{align}
\sup_{x\in U, t\in I} |\f_{t,\e}(x) - \f_t(x) | & \leq C
\e^\g; \label{appr1}\\
\sup_{x\in U, t\in I_0} |\p_t \f_{t,\e}(x)| & \leq
C\e^{\g-1}; \label{appr2}\\
\sup_{x\in U, t\in I_0} |\n_x \f_{t,\e}(x) | & \leq C.\label{appr3}
\end{align}
Let us fix a non-negative function $\eta \in C_0^\infty(\R^{n-k})$
with $\eta = 1$ on $B_{2C}$ and $\eta = 0$ on $\R^{n-k}
\backslash B_{3C}$. We consider the following cut-off function
$$
\chi_\e(t,x) = 1 - \eta\left( \frac{1}{\e^\g}
(\f_{t,\e}(x))_{k+1},...,\frac{1}{\e^\g} (\f_{t,\e}(x))_{n}\right),
$$
for $x\in U$ and $t \in I$. Notice that as $\e \ra 0$ $\chi_\e \ra
1$ for all $t$ and a.e. $x$. Furthermore, due to \eqref{appr1},
$\supp(\chi_\e \phi)$ does not intersect the set $S$ on the time
interval  $I$. Finally, put
$$
\phi_{\e} = \chi_\e \phi.
$$
Due to regularity of $u$ away from $S$, Lemma \ref{L:reg} applied to produce
\begin{multline}\label{localeee}
 \int_{U \times \{t''\}} |u|^2 \phi_\e -  \int_{U \times \{t'\}} |u|^2  \phi_\e - \int_{ U\times [t',t'']} |u|^2 \p_t  \phi_\e
 \\ =  \int_{U\times [t',t'']} (|u|^2 + 2p) u\cdot \n \phi_\e.
\end{multline}
Let us examine the terms in the limit as $\e \ra 0$.  Clearly, the
first two terms on the right hand side will converge to their
natural limits. As to the third term, we have $\p_t \phi_\e = \phi
\p_t \chi_\e$, and
$$
\p_t \chi_\e = - \frac{1}{\e^\g} \sum_{j=k+1}^n \p_t(\f_{t,\e}(x))_j
\p_j \eta\left( \frac{1}{\e^\g}
(\f_{t,\e}(x))_{k+1},...,\frac{1}{\e^\g} (\f_{t,\e}(x))_{n}\right).
$$
Notice that $\p_t \chi_\e$ is supported on the set
$$
U \cap (\f_{t,\e})^{-1}(\R^k \times [-3C\e^\g,3C\e^\g]^{n-k} \cap
B_1),
$$
which is a subset of
$$
A_{\e} = U \cap (\f_{t})^{-1}(\R^k \times [-4C\e^\g,4C\e^\g]^{n-k}
\cap B_1).
$$
We have $|A_{\e}| \sim \e^{(n-k)\g}$. In view of \eqref{E:expon} and \eqref{appr2}
we obtain
\begin{align}
\left|  \int_{U\times [t',t'']} |u|^2 \phi  \p_t \chi_\e \right| &
\leq \frac{|A_{\e}|^{(q-2)/q} }{\e} \int_{t'}^{t''} \left(  \int_{A_{\e}} |u|^q dx
\right)^{2/q} dt \label{notoptimal1}
\\
& \leq \int_{t'}^{t''} \left(  \int_{A_{\e}} |u|^q dx \right)^{2/q}
dt \underset{\e \ra 0}{\longrightarrow} 0.
\end{align}
Let us now examine the right hand side of \eqref{localeee}. We have
$$
u \cdot \n \phi_\e = \phi u \cdot\n \chi_\e + \chi_\e u \cdot \n \phi.
$$
Clearly we can pass to the limit in the integral containing the
second term. As to the first term we have
$$
\n \chi_\e = - \frac{1}{\e^\g} \sum_{j=k+1}^n \n(
\f_{t,\e}(x))_j \p_{j} \eta \left(
\frac{1}{\e^\g}(\f_{t,\e}(x))_{k+1},...,\frac{1}{\e^\g}
(\f_{t,\e}(x))_{n}\right),
$$
which is supported on the set $A_{\e}$. Thus,
\begin{align}
\left| \int_{U\times [t',t'']} (|u|^2 + 2p) u\cdot \n \phi_\e
\right| & \leq \frac{|A_{\e}|^{(q-3)/q}}{\e^\g} \int_{t'}^{t''} \left(  \int_{A_{\e}}
|u|^q dx \right)^{3/q} dt \label{notoptimal2}
 \\
& \sim \int_{t'}^{t''} \left(  \int_{A_{\e}} |u|^q dx \right)^{3/q}
dt \underset{\e \ra 0}{\longrightarrow} 0.
\end{align}
This finishes the proof of \thm{T:main}.
\end{proof}

\subsection{An Onsager-critical improvement}\label{SS:oo}

Let us consider physical units of velocity -- $U$, length -- $L$ and time -- $T$. Then the dimension of the
regularity space $\reg$ is $T^{1/3} U L^{\frac{n-1}{3}}$. We call functional spaces of this dimension Onsager-critical. In the case of point singularities, i.e. $k=0$, \thm{T:main} yields the Onsager-critical condition $u \in L^3 L^{\frac{3n}{n-1}}$
with $\g$ being at least $\frac{3}{n+2}$. Under these circumstances we expect
our result to be optimal. However, this is not the case if $k >0$, since the dimension of $L^3L^{\frac{3(n-k)}{n-k-1}}$ is $T^{1/3}U L^{\frac{n(n-k-1)}{3(n-k)}}$.
Onsager-critical spaces for $k>0$ can be defined using mixed $L^p$ spaces relative to the slices $S(t)$. Assuming that each $S(t)$ is a $k$-dimensional
smooth submanifold of $\R^n$ we consider local normal fiber bundle $S^\perp(t)$. Thus, each fiber $S^\perp(x,t)$ is a $\gamma$-smooth in time local tile orthogonal to the surface $S(t)$. We can now define the local space $u \in L^3_t L^p_{S} L^q_{S^\perp}$ by requiring over coordinate neighborhood $U \times I$ the condition
$$
\int_I \left(  \int_{S(t) \cap U} \left(   \int_{S^\perp(x,t) \cap U} |u(x,y,t)|^q d\s^{n-k}_t(y)  \right)^{p/q} d \s^k_t(x)   \right)^{3/p} dt  < \infty,
$$
where $d\s_t$ indicates the surface measure of the corresponding dimension. Notice that the space $L^3_t L^3_{S} L^{\frac{3(n-k)}{n-k-1}}_{S^\perp}$
is in fact Onsager-critical. In general, Theorem \ref{T:main} can be restated by requiring
\begin{equation}\label{optim}
u \in L^3_t L^3_{S} L^q_{S^\perp}
\end{equation}
under the same assumptions on $n,k,\g,q$. In particular, we obtain energy conservation if
\begin{equation}\label{optim2}
u \in L^3_t L^3_{S} L^{\frac{3(n-k)}{n-k-1}}_{S^\perp} \text{ and } \gamma \geq \frac{3}{n-k+2}.
\end{equation}

In order to reprove \thm{T:main} under new condition \eqref{optim} one has to
simply apply the H\"{o}lder inequality in \eqref{notoptimal1} and in \eqref{notoptimal2} only to the integrals over $S^\perp(x,t)$, the rest of the argument being the same. We leave details for the reader.

\subsection{Other extensions}\label{SS:ext}

Since our argument is local, it is readily extendable to the case of locally finite union of singular sets. Specifically, suppose that in every coordinate neighborhood $V = U \times I$
\begin{equation}\label{sun}
S = \bigcup_{j = 1}^{N_V} S_j,
\end{equation}
where $S_j$'s are $k_j$-dimensionally $C^{\g_j,1}$-covered in $V$. We can use the product of cut-offs
$$
\chi_\e = \prod_{j= 1}^{N_V} \chi_\e^j
$$
to run the argument. The conclusions of \thm{T:main} remains true under the corresponding assumptions on $u$ locally near each $S_j$. The result of Section \ref{SS:oo} can be modified similarly.

\def \norm {\vec{\nu}(x,t)}
\def \vnu  {\vec{\nu}}

\section{The case of hypersurface: slits} \label{S:slits}

In this section we will study the case $k=n-1$. We will assume
special geometric properties of the singular set $S$. Namely, let
$S$ be a $C^1$-family of closed orientable $C^2$-submanifolds of $\R^n$. For every
$(x_0,t_0) \in S$ there exist $U$, $I$ and a local parametrization
$r = r(\bar{y},t)$ of $S(t) \cap U$ for all $t \in I$, where $r \in C^{2,1}_{\bar{y},t}$,
and $\bar{y} = (y_1,...,y_{n-1}) \in B_1^{n-1}$. Let $\norm$ be the positively oriented unit normal to $S(t)$.
We consider a coordinate system on a smaller neighborhood that is most suitable for dealing with
normal limits. For $\e_0
>0$ small we define
$$
\psi_t(\bar{y},y_n) = r(\bar{y},t) + \e_0 y_n \vnu(r(\bar{y},t),t),
$$
for $|y_n| < 1$. Since $S$ is sufficiently smooth, this defines a
diffeomorphism of $B_1^n$ onto an open neighborhood $U(t)$ with $S(t)
\cap U(t) = S(t)\cap U$ for all $t\in I$. It will be convenient in the future to deal with
$U$ independent of $t$. So, reducing the time interval
if necessary we can find a new neighborhood $U \ss U(t)$ for all
$t\in I$, such that
$$
\psi_t( (-1,1)^{n-1} \times (-\e_1,\e_1) ) \ss U \ss \psi_t( (-1,1)^{n-1} \times (-\e_2,\e_2) )
$$
for all $t\in I$ and some $c_2>c_1>0$. The direct product $V = U\times I$ along with the map $\f_t = \psi_t^{-1}$
define a new coordinate chart containing $(x_0,t_0)$. Let us also define
the normal segments for every $(x,t) \in V$:
\begin{align*}
\G_+(x,t) & = (x,t) + \vnu(x,t)[0,\e_1] \\
\G_-(x,t) & = (x,t) + \vnu(x,t)[ - \e_1, 0].
\end{align*}
We may further truncate the segments to ensure that for some open neighborhood $W$ of $S$ we have
$\bigcup_{S}\G_\pm \ss W$. For a function or field $f$ on $W$ we denote by $f^*_\pm:S \ra \R$
the normal maximal function defined by
$$
f^*_\pm(x,t) = \sup_{x' \in \G_\pm(x,t)} |f(x',t)|,
$$
and by $f_\pm$ the limits
$$
f_\pm(x,t) = \lim_{x' \ra x, x'\in \G_\pm(x,t)} f(x',t),
$$
if the latter exist.

We now introduce a measure on each $S(t)$ whose role will be clear
in a moment. We start by defining it locally on every chart $U\cap
S(t)$. For this purpose let us fix a scalar-valued function $H(x,t)
\in C^1$ with level surface $\{H(x,t) = 0\} = S(t) \cap U$ for all
$t\in I$, and such that $\n_x H \neq 0$ agrees
 with $\vnu$. For instance, $H = (\f_t(x))_n$. Let us consider the measure
\begin{equation}\label{meas}
d \mu_t^{U}(x) = \frac{\p_t H}{|\n_x H|} d \s_t(x)
\end{equation}
where $d \s_t(x)$ is the surface measure of $S(t)$.
Notice the following identities
\begin{align*}
H(r(\bar{y},t),t) & = 0 \\
\p_t H(x,t) & = - \p_t r(\bar{y},t) \cdot \n_x H(x,t),
\end{align*}
where $x = r(\bar{y},t) \in S(t)\cap U$. Thus, in local coordinates,
$$
d \mu_t^{U}(x) = - \p_t r(\bar{y},t) \cdot
\vec{\nu}(r(\bar{y},t),t) J_t(\bar{y}) d\bar{y},
$$
where $J_t(\bar{y})$ is the volume element. We see that the
definition of $d \mu_t^{U}$ is independent of $H$. Yet \eqref{meas} shows that it is also independent of particular parametrization of $S(t)$.
Now, let $f \in C_0(S(t))$ be a continuous function with compact support on
$S(t)$. Arguing as in the proof of \thm{T:main} we find a finite
cover of $\supp(f)$ by $\{ U_i\}_{i=1}^M$ with the corresponding
 partition of unity $\{ g_i\}_{i=1}^M$  over  $\supp(f)$. Define
\begin{equation}\label{globmeas}
\int_{S(t)} f d \mu_t(x) = \sum_{i=1}^M \int_{S(t)} f g_i d
\mu_t^{U_i}(x).
\end{equation}
This is a well-defined measure over $S(t)$. For instance, if $S(t)$
is given by the graph of a periodic in spacial variables function
$x_n = z(x_1,...,x_{n-1},t)$, then
$$
d \mu_t = - \p_t z(x_1,...,x_{n-1},t) dx_1...dx_{n-1}.
$$

The measure $d \mu_t$ arises naturally in the following calculation.
Let us fix a coordinate chart $(V,\f_t)$ as above, define $\eta$ as in
the previous section with $k = n-1$, and denote 
\begin{equation}\label{chiagain}
\chi_\e(x,t) = 1- \eta(\e^{-1} ( \f_t(x))_n),
\end{equation} 
for $\e < \e_1$.
\begin{lemma}\label{L:microlim}
Let $f:V \ra \R$ and $u: V \ra \R^n$ be such that the limits $f_\pm(x,t)$ and
$u_\pm(x,t)$ exist for a.e. $t\in I$ and a.e. $x\in S(t)$ with
respect to $d\s_t$, and $f^*_\pm,\,u^*_\pm \in L^1(d\s_t dt)$. Then
\begin{align}
\lim_{\e \ra 0} \int_{V_0} f \p_t \chi_\e dx dt & =
\int_I \int_{S(\t)}(f_+ - f_-) d\mu_\t d\t, \label{microlim}\\
\intertext{and} \lim_{\e \ra 0} \int_{V_0} u \cdot \n_x \chi_\e dx
dt & = \int_I \int_{S(\t)} (u_+ - u_-)\cdot \vnu d\s_\t d\t.
\label{microlim2}
\end{align}
\end{lemma}
\begin{proof}
Let us denote $H(x,t) = (\f_t(x))_n$.
To prove \eqref{microlim} let us observe
$$
\int_{V} f \p_t \chi_\e dx dt = - \int_V \e^{-1} \eta'(\e^{-1} H(x,t)) \p_t H(x,t) f(x,t) dxdt.
$$
As a guiding point we recall the classical microlocal limit
$$
\frac{1}{\e} \int_{ 0 \leq H \leq \e} g dx \ra \int_{H = 0} \frac{g}{|\n H|} d\s.
$$
By changing the variables we obtain the integral
$$
\int_{V} f \p_t \chi_\e dx dt =- \int_{I\times B^{n-1}_1} F_\e(\bar{y},t) d \bar{y} dt,
$$
where
\begin{multline*}
F_\e(\bar{y},t) = -\int_{|y_n| < \e_1} f(\psi_t(\bar{y},y_n),t) \p_t
H(\psi_t(\bar{y},y_n),t) \e^{-1} \eta'(\e^{-1} y_n) \\ \O_t(\bar{y},y_n)
dy_n dt,
\end{multline*}
and
$$
\O_t(y) = \left| \det{\frac{D\psi_t}{Dy}} \right|.
$$
Given our choice of $H$ we have
\begin{align}
H(\psi_t(\bar{y},y_n),t) & = y_n, \\
\p_t H + \n_xH \cdot \p_t \psi_t(\bar{y},y_n)& = 0.
\end{align}
So, as $y_n \ra 0$ we obtain
$$
\p_t H \ra - \p_t r(\bar{y},t) \cdot \n_x H(r(\bar{y},t))
$$
uniformly in $\bar{y} \in B^{n-1}_1$. Moreover,
$$
\O_t(y) \ra  \e_0 J_t(\bar{y}).
$$
Using that $\n_x H(x,t) = \e_0^{-1} \norm$ we obtain the uniform convergence
$$
\p_t H \O_t \ra - \p_t r \cdot \vnu J_t(\bar{y}).
$$
Let us observe now that as $\e$ gets sufficiently small, we have $\psi_t(\bar{y},y_n) \in
\G_{\sgn{y_n}}(\psi_t(\bar{y},0),t)$ for all $\bar{y} \in
B_1^{n-1}$, and $\psi_t(\bar{y},y_n)$ approaches the surface
orthogonally. The condition $f^*_\pm \in L^1(d\s_t dt)$ implies that all $F_\e$ have a common integrable majorant.
This enables us to pass to the limit and arrive at \eqref{microlim}.
The proof of \eqref{microlim2} is similar.
\end{proof}

\begin{definition}\label{D:slit}
Let $u$ be a weak solution to the Euler equations. The surface $S$ is called a \emph{slit}
of $u$ if
\begin{itemize}
\item[1)] The limits $u_\pm, p_\pm$ exist for a.e. $t \in [0,T]$ and a.e. $x\in S(t)$,
\item[2)] $u^*_\pm \in L^2(d\s_t dt)_{\loc}$ and $p^*_\pm \in L^1(d\s_t dt)_{\loc}$
\end{itemize}
\end{definition}
\begin{lemma}\label{L:slit}
Let $u$ be a weak solution to the Euler equations, and $S$ be a slit. Then the following is true:
\begin{itemize}
\item[1)] $u_+\cdot \vnu = u_-\cdot \vnu:=u_\nu$ and $p_+ = p_-$ for a.e. $t\in [0,T]$ and a.e. $x \in S(t)$;
\item[2)] $d\mu_t + u_\nu d\s_t = 0$ for a.e. $t$ on the set $u_+ \neq u_-$.
\end{itemize}
\end{lemma}

\begin{proof}
As before we reduce the statements of the lemma to the local
coordinate neighborhood $V = U \times I$ defined earlier. Let us consider an
arbitrary scalar function $g\in C_0^\infty(V)$. From the
divergence-free condition on $u$ we obtain
$$
\int_{V} u \cdot \n(g \chi_\e) = 0.
$$
Letting $\e \ra 0$ we obtain from Lemma \ref{L:microlim}
$$
\int_{V} u \cdot \n g + \int_I \int_{S(\t)} g( u_+\cdot \vec{\nu}
- u_- \cdot \vec{\nu}) d\s_\t d\t = 0.
$$
Using the divergence-free condition again and the free choice of $g$
we obtain
\begin{equation}\label{normal}
u_+\cdot \vnu = u_-\cdot \vnu.
\end{equation}

Consider an arbitrary vector-valued function $a \in C_0^1(V)$, and
$\psi = a \chi_\e$. By continuity, the regularity of $\psi$ is
sufficient to substitute $\psi$ into \eqref{weaksolalt}. We obtain the
following identity:
\begin{multline*}
     - \int_{V} u\cdot
    \p_\t a \,\chi_\e - \int_{V} u\cdot a
    \p_\t \chi_\e = \int_{V} (u\otimes u) : \n a\, \chi_\e + \int_{V} (u\cdot a) (u\cdot \n \chi_\e)\\ + \int_{V} p(a \cdot \n \chi_\e + \chi_\e \diver a).
     \end{multline*}
    Using \eqref{microlim} and \eqref{microlim2} we obtain in the limit as $\e \ra 0$
\begin{multline*}
     - \int_{V} u\cdot
    \p_\t a  - \int_I \int_{S(\t)} (u_+ - u_-) \cdot a
    d\mu_\t d\t \\= \int_{V} (u\otimes u) : \n a + \int_I \int_{S(\t)} (u_+ - u_-)\cdot a  u_\nu  d\s_\t d\t \\+
    \int_I \int_{S(\t)} (p_+ - p_-) a_\nu d\s_\t d\t +   \int_{V} p \diver a.
     \end{multline*}
Using the identity for the weak solutions \eqref{weaksolalt} with $\psi = a$ we see that only the boundary terms remain:
\begin{multline*}
    - \int_I \int_{S(\t)} (u_+ - u_-) \cdot a
    d\mu_\t d\t = \int_I \int_{S(\t)} (u_+ - u_-)\cdot a  u_\nu  d\s_\t d\t \\+ \int_I \int_{S(\t)} (p_+ - p_-) a_\nu d\s_\t d\t .
 \end{multline*}
 Let us choose $a$ of the form $a = \vec{\nu} g$, where $g\in C_0^1(V)$. Using \eqref{normal} we have
 $$
 \int_I \int_{S(\t)} (p_+ - p_-) g d\s_\t d\t = 0.
 $$
 This readily implies $p_+ = p_-$ a.e. Going back to the previous identity we notice that 2) holds as well due to arbitrariness of $g$.
\end{proof}

\begin{theorem}\label{T:slit}
Suppose that $u\in L^3(\R^n\times [0,T])$ is a weak solution to the Euler equations and the singular set $S$ of $u$ is a slit.
Suppose further that $u^*_\pm \in L^3(d\s_t dt)_\loc$,
$p^*_\pm \in L^{3/2}(d\s_t dt)_\loc$. Then $u$
conserves energy.
\end{theorem}

In view of our discussion in Section \ref{SS:oo} we notice that the conditions of \thm{T:slit} are Onsager-critical. We therefore
expect these conditions to be optimal as far as our argument in concerned.

\begin{proof}
As in the proof of \thm{T:main} we reduce the problem to proving the
local energy equality \eqref{localee}. As before  $H(x,t) = (\f_t(x))_n$ and $\phi_\e$ is defined by \eqref{chiagain}. The regularity of $u$ away from the slit $S$ enables us to use Lemma \ref{L:reg} with $\phi_\e$.
Using the results of Lemmas \ref{L:microlim}
and \ref{L:slit} we can pass to the limit as $\e \ra 0$ and obtain
\begin{align*}
\int_{U \times [t',t'']} |u|^2 \p_t \phi_\e & \ra
\int_{U \times [t',t'']} |u|^2 \phi \\& + \int_{t'}^{t''}
\int_{S(\t)}
(|u_+|^2 - |u_-|^2) \phi d \mu_\t d\t, \\
\intertext{and} \\
\int_{U \times [t',t'']} (|u|^2 + 2p) u \cdot \n \phi_\e & \ra
\int_{U \times [t',t'']} (|u|^2 + 2p) u\cdot \n \phi \\&+
\int_{t'}^{t''} \int_{S(\t)} (|u_+|^2 - |u_-|^2) u_\nu d \s_\t d\t
\\ &+ 2 \int_{t'}^{t''} \int_{S(\t)} (p_+ - p_-)u_\nu d \s_\t d\t .
\end{align*}
According to Lemma \ref{L:slit} the surface integral terms sum up to
zero, and \eqref{localee} follows.
\end{proof}

Arguing as in Section \ref{SS:ext} we can include the result of \thm{T:slit} in obtaining more general singular set configurations.
Thus, the union \eqref{sun} may involve finitely many slits accompanied by the corresponding conditions on $u$ and $p$.

We remark that one can also state the conditions of \thm{T:slit} and Definition \ref{D:slit} in terms of more conventional 
non-tangential limits and maximal functions. It would be interesting to know whether the condition $u^*_\pm \in L^{2q}(S)$ automatically implies $p^*_\pm \in L^q(S)$.

\section{Energy of vortex sheets}\label{S:sheets}
Naturally, the conditions of \thm{T:slit} apply to vortex sheet solutions.
Vortex sheets in the classical sense (as opposed to those defined by Delort \cite{Delort}) are singular solutions to the Euler equations
with vorticity concentrated on a hypersurface (see \cite{Saffman}). For notational convenience we will consider the two dimensional case, although
all what follows holds true in three dimensions as well. In 2D a vortex sheet is described by the graph of a regular function $\zeta(\a,t) = (\a, h(\a,t))$ and vorticity density $\g = \g(\a,t)$ on the graph. Typically, one assumes $2\pi$-periodicity on $h$ and $\g$. Thus, in complex variable notation the velocity field off the sheet is given by the Biot-Savart law
$$
\bar{u}(z,t) = \frac{1}{4\pi i} \int_{-\pi}^\pi \cot\left(\frac{z-\zeta(\a,t)}{2} \right) \g(\a,t) d \a.
$$
Provided $\g$ has enough smoothness on a time interval $[0,T]$, the standard potential theoretical considerations imply that $u\in L^\infty_t L^\infty_x$,
the non-tangential, and hence normal, limits exist and are given by
$$
u_\pm(\a,t) = - \frac{1}{4\pi i} PV \int_{-\pi}^\pi \overline{\cot}\left(\frac{\zeta(\a,t)-\zeta(\a',t)}{2} \right) \g(\a',t) d \a' \mp \g(\a,t)\vec{s},
$$
where $\vec{s}$ is the unit tangent vector oriented in the positive direction of the $x$-axis. The pressure can be recovered from Bernoulli's function, and is given by the double-layer potential formula
$$
p = - \frac{1}{2} |u|^2 + \frac{1}{2} \mathcal{D}( |u_+|^2 - |u_-|^2  ).
$$
From the classical jump relations for the double-layer potential $\mathcal{D}$ we conclude that the limits $p_\pm$ exist,
$p_+ = p_- = \frac{1}{4}(|u_+|^2 +|u_-|^2 )$ and $p^*_\pm \in L^q(d\s_t dt)_\loc$ for all $1 \leq q<\infty$. Thus, according to Definition \ref{D:slit}
the classical vortex sheet is a slit. The equation 2) in Lemma \ref{L:slit} is nothing but the well-known evolution law of the sheet:
$$
\p_t h = -U_1 \p_\a h + U_2,
$$
where $U = \frac{1}{2}(u_+ + u_-)$.
In order for the total kinetic energy of the vortex sheet to be finite we assume vanishing of the total circulation:
$$
 \int_{-\pi}^\pi \g(\a,t) d\a = 0.
$$
Under this condition, $u\in L^\infty L^2$. By interpolation with  $u\in L^\infty L^\infty$ we obtain
$u\in L^3L^3$. Therefore, the conditions of \thm{T:slit} are satisfied and we arrive at the following corollary.
\begin{corollary}\label{C:sheet}
Suppose that $\g,h \in C^\infty([0,T] \times [-\pi,\pi])$, and the total circulation of $\g$ is zero. Then the energy of the vortex sheet is conserved.
\end{corollary}
Vortex sheets of this nature are known to exist in 2D and 3D locally in time in spaces of functions that admit analytic extension to a complex strip (see \cite{Caflisch,SSBF}). In general, the global existence is precluded by occurrence of the roll-up singularity (see \cite{Krasny}). The conditions on Cauchy data stated in  \cite{SSBF}  that guarantee local existence allow for sheets with zero circulation. Thus, Corollary \ref{C:sheet}
applies to a variety of existing vortex sheets. However, the proof of \thm{T:slit} applies to obtain local energy balance relation for sheets with infinite energy as well.

%\bibliographystyle{plain}
%\bibliography{ee4Euler}

\end{document}